\documentclass[12pt]{article}

\usepackage[english]{babel}
\usepackage[utf8]{inputenc}

\usepackage[nofiglist,notablist]{endfloat}

\usepackage[intlimits]{amsmath} 
\usepackage{amsthm,mathrsfs}    
\usepackage[comma,authoryear]{natbib}
\usepackage{url}
\usepackage{color}
\numberwithin{equation}{section}
\usepackage{amssymb,amsmath}
\usepackage{subfigure}
\usepackage{bm}
\usepackage{booktabs} 
\usepackage{graphicx,graphics,epsfig}
\usepackage{float}
\usepackage{multirow}

\theoremstyle{plain}
\newtheorem{thm}{Theorem}

\newcommand{\bthm}{\begin{thm}}
\newcommand{\ethm}{\end{thm}}

\newcommand{\bpf}{\begin{proof}}
\newcommand{\epf}{\end{proof}}

\newcommand{\taubh}{\tau_{{\rm BH}}}

\theoremstyle{definition}

\usepackage[a4paper,margin=.9in,centering,text={16.2cm,24cm}]{geometry}

\usepackage[compact]{titlesec}
\newcommand{\bib}{\bibliography{ref-bib}\bibliographystyle{biom}}
\usepackage{deep-macro}

\begin{document}
\begin{center}
{\Large {\bf Large Scale Signal Detection: A Unified Perspective}}
\\[.2in]
Subhadeep Mukhopadhyay\\
Department of Statistical Science, Temple University \\ Philadelphia, Pennsylvania, 19122, U.S.A.\\
\texttt{deep@temple.edu}\\[.65em]
\emph{Dedicated to John Tukey, on the occasion of his 100th birthday.} \vskip.65em
{\bf ABSTRACT}\\
\end{center}
\vspace{-.7em}
There is an overwhelmingly large literature and algorithms already available on `large scale inference problems' based on different modeling techniques and cultures. Our primary goal in this paper is \emph{not to add one more new methodology} to the existing toolbox but instead (a) to clarify the mystery how these different simultaneous inference methods are \emph{connected}, (b) to provide an alternative more intuitive derivation of the formulas that leads to \emph{simpler} expressions in order (c) to develop a \emph{unified} algorithm for practitioners. A detailed discussion on representation, estimation, inference, and model selection is given. Applications to a variety of real and simulated datasets show promise. We end with several future research directions.

\vspace*{.15in}
\section{Introduction}
\subsection{Mixed Sample Inference}
Consider $n$ i.i.d. observation $Z_1,\ldots,Z_N$ from the continuous distribution $F=\pi_0 F_0 + (1-\pi_0) G$, $0<\pi_0 \le 1$, where $F_0$ is specified distribution (known) and $G$ is unknown mixing distribution. As the given random sample may contain observations from both $F_0$ and $G$, we call this framework a `mixed-sample problem.' This can be considered as an intermediate statistical inference problem between two extremes: one and two sample inference problems. Typical goals of mixed data inference problems include: (1) identifying the $Z_i$'s coming from the distribution $G$ (also known as signal detection), and (2) estimating the proportion $\pi_0$.
\vskip.5em

\emph{Motivation}. A multiple hypothesis testing problem can be thought of as a mixed-sample problem by considering $Z_i$'s to be the test statistic corresponding to a null hypothesis $H_{0i}$ ($i=1,\ldots,N$), where the goal is simultaneous inference. Our main motivation is to develop a comprehensive modeling framework for big-$N$ signal detection (or large-scale mixed-sample) problems in a way that is \emph{analogous} to the one-sample and two-sample approaches so that we can cover the whole statistical inference spectrum (one sample $\rightarrow$ mixed sample $\rightarrow$ two sample) using one general concept and tool -- which will simplify theory, practice and teaching.
\subsection{Goals and Organization}
The topic of multiple hypothesis testing, which began with the pioneering work of Tukey on ``Comparing individual means in the analysis of variance'' published in $1949$ Biometrics\nocite{tukey1949mc}, gained new momentum in the $21$st century with modern high-throughput data acquisition. Since then an enormous amount of research has been conducted on this topic.  Following is a summary of our research contributions:

(1) We seek to \emph{condense and systemize} the vast literature on large-scale simultaneous inference methods based on a few fundamental concepts and notations introduced in Section 2. The author believes that this could ultimately help applied data scientists to gain better insight for selecting the proper signal detection tool from the existing large inventory by simplifying the practice. Section 2.2 describes a \emph{unified representation theory and an alternative rationale} of three main classes of large-scale inference procedures: Benjamini Hochberg's false discovery rate (FDR) approach \citep{BH95},  Efron's empirical Bayes local FDR proposal \citep{efron01}, and Donoho Jin's Higher Criticism thresholding \citep{donoho2004}.

(2) Section 2.2 further shows how this new theoretical framework for large scale simultaneous inference problems can \emph{integrate and reconcile} frequentist and Bayesian cultures, \emph{bypassing} the old philosophical and ideological differences \citep{benj08} to develop unified efficient algorithms using single concept and notation.

(3) The theory presented in Section 2 motivated us to introduce a \emph{new class of nonparametric signal detection algorithm} by converting the large-scale inference problems into a function estimation problem based on the ideas of comparison density and distribution (CD); see Web Appendix A for more discussion on the historical significance of the current research. Section 3.1 describes the statistical modeling challenges. A comprehensive CD-based nonparametric modeling solution for `signal-hunting' procedure is described in Section 3.2-3.5, which consists of three main components: (A)
We discover that comparison density takes a `\emph{universal stylized shape}' for almost all large-scale inference problems - `U' shaped density over the compact unit interval. To tackle this nonstandard density we propose a specially designed Skew-Beta decomposition-based nonparametric estimation algorithm; (B) To model empirical null we propose a robust biweight M-estimator, which is easy to compute; (C) Finally, we describe the Minimum Deviance Criteria (MDC) algorithm for estimating true null proportion. Section 4 presents numerical results, and discussions are given in Section 5.

%
\section{Theory}
\subsection{Background on Functional Statistical Inference}
Let $Z_1,Z_2,\ldots,Z_N$ be a mixed random sample, with the majority of the observations coming from continuous null distribution $F_0$, and a small proportion from unknown contaminating distribution $G$: $F=\pi_0F_0 + (1-\pi_0) G,\, 0< \pi_0 \le 1$, where $Z_i$'s are the corresponding test statistic for the hypothesis testing problem $H_{0i}$ and the goal is simultaneous inference. Here we will provide an intuitive introduction to the functional statistical approach for large-scale simultaneous inference problems whose goal is to detect false null hypotheses.

The cornerstone of our approach is based on the concept of comparison density \citep{parzen83a}. For continuous $F$ and $G$, define comparison distribution function $D(u;G,F):=F(G^{-1}(u))$ and the corresponding comparison density by:
\beq \label{eq:cd}
d(u;G,F)\,=\, \dfrac{f(G^{-1}(u))}{g(G^{-1}(u))}, \quad 0 < u<1.
\eeq

Under $H_0:F=F_0$, i.e., when all the observations actually come from $F_0$, the theory of weak convergence of comparison distribution process \citep{parzen99,Thas10} tells us $\sqrt{N}\big( \widetilde D(u) - u\big)  \xrightarrow{d} \mathbb{B}(u)$ as $N \rightarrow \infty$, where $\mathbb{B}(u) ,0<u<1$ is the Brownian bridge process, and $\widetilde D(u)= \widetilde F(F_0^{-1}(u))$, where $\widetilde{F}(z;Z)=N^{-1}\sum_{i=1}^N \mathbb{I} (Z_i \le z)$. In mixed-sample problems, our goal is to identify the signals or the $Z_i$'s that are coming from $G$. If all the null-hypothesis $H_{01},H_{02},\ldots, H_{0N}$ are true (i.e., under $H_0:F=F_0$), we would expect $\widetilde D(u):= D(u;F_0,\widetilde F) \approx u$. Thus, intuitively, the collection of $u$'s, for which the distance between $\widetilde D(u)$ and $u$ are large, should be suspected as potential signals. In other words, \[\{u:\widetilde D(u)/u\,>\,\ga\}, \] for a suitably selected threshold might contain the signals or the rejected null hypotheses. Now instead of comparison distribution function $D(u)$, in the same spirit, we can develop a similar strategy for detecting false null hypotheses based on comparison density $d(u)$. So the performance of simultaneous inference critically depends on the fundamental comparison density (or equivalently distribution) function.

\emph{From Multiple Testing to Nonparametric Function Estimation}. The crux of this section is this: the notion of comparison distribution and density allows us to transform the simultaneous hypothesis testing problem into a nonparametric function approximation problem - a functional statistical approach to large scale inference. In the next section, we will make this intuitive  mixed-sample signal detection approach formal and rigorous. 
\subsection{Unified Representation and Theory of Threshold Selection}
The previous section intuitively argued how to build a CD-based simultaneous inference algorithm. The purpose of this section is threefold: (1) to derive proper thresholds for CD-based simultaneous inference procedures that can lead to desirable methods with provable guarantees; (2) to show how the Frequentist and Bayesian large-scale inference algorithms can be connected using the theory of reproducing kernel Hilbert space (RKHS) of the limiting Brownian bridge process of the comparison distribution; (3) to simplify the currently existing vast literature on large-scale multiple testing by providing unified representation based on a single concept and notation -- comparison density/distribution, to our knowledge largely an unexplored area of research.

Let's start with the question of how to select the proper threshold. The answer depends on the appropriate error rate that we seek to control in multiple hypothesis testing. Based on particular set up and purpose of any specific inference problems \cite{benjamini2010} provided a list of $14$ such error rates.

One of the important error-rate criteria (among many possibilities) in large-$N$ multiple testing is to control the false discovery rate (FDR) -- the expected proportion of Type I errors among the rejected hypotheses, suggested by \cite*{BH95}. The following theorem describes how to select the threshold for comparison distribution-based method that can guarantee FDR $\le \al$ (for some preselected significance level $\al$) under any arbitrary $\pi_0$ and $G$. A sketch of the proof is provided in Web Appendix B.

\begin{thm}
Consider testing $N$ independent null hypothesis $H_{(01)},\ldots,H_{(0N)}$ based on corresponding ordered p-values $u_{(1)},\ldots,u_{(N)}$. Define the index set
\beq \label{eq:CDBH}
\mathcal{R} =\Big\{i \le k: k= \argmax_{i} \dfrac{\widetilde{D}(u_{(i)})}{u_{(i)}}\,\ge\, \dfrac{\eta}{\al}\Big\}.
\eeq
Then the procedure that rejects $H_{(0i)}, i \in \mathcal{R}$ controls FDR at the level $\al$, regardless of the distribution of the test statistic corresponds to false null hypothesis.
\end{thm}

To see how CD-based algorithm \eqref{eq:CDBH} essentially equivalent to the (adaptive) Benjamini Hochberg's (BH) FDR controlling procedure, note that $\widetilde{D}(u_i;F_0,\widetilde F)= \widetilde F F_0^{-1}(u_{(i)})=i/N$, which is equivalent to saying that reject $H_{(0i)}$ $i=1,\ldots,k$, where
\[k~=~\max\big\{1 \le i \le N: u_{(i)} \le \dfrac{i}{N} \dfrac{\al}{\pi_0} \big\}.\]

\emph{RKHS-based dual representation}. Here we will give an alternative CD-based FDR controlling procedure motivated by the isometric isomorphism  (one-one inner product preserving map) relation between the Hilbert space spanned by the Brownian bridge process (which is the limit process of $\sqrt{N}\big( \widetilde D(u) - u\big)$ and the corresponding RKHS associated with the covariance kernel $K_{\mathbb{B}}(u,v)$. The following theorem provides the complete characterization of the RKHS.

\begin{thm}
The reproducing kernel Hilbert space $H_K$ associated with the Brownian bridge covariance kernel $K_{\mathbb{B}}(u,v)= \min(u,v)-uv,\, 0<u,v<1$ consists of $L^2$ differentiable functions with the inner produce $\big\langle \phi,\psi \big\rangle=\int_0^1 \phi'(u) \psi'(u) \dd u$ that satisfies $\phi(0)=\phi(1)=0$.
\end{thm}

It is sufficient to prove the reproducing property relative to the inner product in the following sense, for all $\phi \in H_K$
\beq \big\langle K_{\mathbb{B}}(u,\cdot), \phi\big\rangle\,=\,\phi(u), \text{~for any\,} u \in [0,1]. \eeq
For details of the proof see Web Appendix C.
\vskip.5em

\emph{Moving toward Empirical Bayes Formulation}. As the RKHS associated with the Brownian bridge covariance kernel equipped with the norm squared $\|h\|^2=\int_0^1|h'(u)|^2 \dd u$,  this tells us that the other bonafide signal detection technique can equivalently be constructed by looking at the distance between $D'(u)=d(u)$ and $1$. That is, the collection of $u$'s for which $d(u;F_0,F)$ \textit{substantially} deviates from $\rm{Uniform}[0,1]$, are precisely the non-null candidates that we are searching. A comparison density $d(u;F_0,\widetilde F)$  based signal characterization  (instead of comparison distribution function $D(u;F_0,\widetilde F)$ is given in the following result, which turns out to be an alternative way of expressing Efron's empirical Bayes local fdr \citep{efron01} formula.
\begin{thm} Local false discovery rate can alternatively be represented using the comparison density
\beq \label{eq:locfdr2} \fdr(z) ~:= ~\Pr\{ {\rm null} \mid Z=z\}~ =~  \dfrac{\pi_0}{d( F_0(z); F_0, F )}.\eeq
\end{thm}
The local fdr is defined as the conditional probability of a case being null or noise given $Z=z$,
\beq \label{eq:locfdr}
\fdr(z)\,=\, \Pr({\rm null} \mid Z=z)\,=\, \pi_0 \,\dfrac{f_0(z)}{f(z)}\,=\, \pi_0/d(F_0(z);F_0,F).
\eeq
It is known that the local FDR is more conservative than the BH procedure. To address this issue, \cite{efron07} recommended the threshold $\{u:d(u;F_0,\widetilde F)>\taubh/2\}$ (under the Lehmann alternatives), where $\taubh=\pi_0/\al$, which he argued controls size, or
Type I errors at the desired level. The goal of local FDR algorithm is to estimate the conditional probability $\Pr({\rm null} \mid Z=z)$ from the data. The traditional procedure estimates \textit{separately} $\hpi_0$, $\hf_0(t)$ and $\hf(t)$, while our comparison density-based approach \emph{directly} estimates the probability; see Web Appendix D for more discussion on this point.

\emph{Relation with Higher Criticism.} Another related method for detecting rare/weak signals in large-scale experiments is ``Higher Criticism'' (HC) thresholding (which is not primarily designed to protect against false discoveries)  introduced by \cite{donoho2004} generalizing Tukey's proposal \citep{tukey89}. Following is the equivalent representation of HC-procedure using our notation. Reject $H_{(0i)}$ $i=1,\ldots,k$ where
\beq \label{eq:HC}
k ~=~\Big\{ 1\le i \le \al_0 N:\, \argmax_{i} \dfrac{\widetilde D(u_{(i)}) - u_{(i)} }{\sqrt{u_{(i)} (1-u_{(i)})}}  \Big\} ,
\eeq
where $\al_0=.1$ is a common choice. We can study the asymptotic properties of HC procedure using the limiting unit variance Gaussian process of $\sqrt{w(u)}\big( \widetilde D(u) - u\big)$, $w(u)=1/u(1-u)$ with the corresponding covariance kernel $K(u,v)=\dfrac{\min(u,v)-uv}{\sqrt{u(1-u)v(1-v)}}$; See \cite{donoho2004} for further details.

\emph{A New Class of Nonparametric Large-scale Inference Procedure.} The representation results proved in this section allow us to develop a new class of CD-based approaches by specifying different models and estimation strategies for $\widehat d$ (or equivalently $\widehat D(u)=\int_0^u \widehat d(v)\dd v$). For example, the \emph{raw-empirical} comparison distribution function $\widetilde D$ generates the BH procedure \eqref{eq:CDBH} or the HC procedure \eqref{eq:HC}, which can be upgraded to a more enhanced technique by \emph{plugging in smooth nonparametric estimate} $\widehat D(u)=\int_0^u \widehat d(v)\dd v$.

\emph{Nonparametric Function Approximation.}
We argue comparison density is an indispensable tool that provides the required theoretical foundation to develop a single general algorithm for modern large-scale signal detection problems. The CD-based approach unifies the frequentist BH procedure and empirical Bayes local false discovery concepts by converting the simultaneous inference problem into a nonparametric function estimation problem.

Any comprehensive multiple testing algorithm should have three main components: (1) $\pi_0$ (null-proportion), (2) $F_0$ (empirical null distribution), and (3) $d(u;F_0,F)$ (comparison density). The efficiency of the signal detection procedure directly depends on the quality of estimation of these quantities. In the next section, we will describe a nonparametric estimation algorithm (CDfdr) in a step-by-step manner (along with real data illustration) to achieve this goal.

The following section starts with an interesting observation that virtually all large-scale inference problems generate comparison density with a very \emph{special universal shape} (see Web Appendix K for more discussion on the tail-behaviour of mixed-sample comparison density). It turns out that the traditional nonparametric density estimation methods fail to satisfactorily capture this typical shape of the comparison density over a compact $[0,1]$ support. For that purpose, we develop a \emph{specially designed nonparametric algorithm} (tailored for this particular characteristic shape) for parsimonious comparison density estimation $\widehat{d}$, which will be used for detecting signals or false null hypotheses and to estimate the null proportion $\pi_0$.

\section{Estimation}
\subsection{Two Modeling Challenges}
\begin{figure*}[!ht]
 \centering
\includegraphics[height=.45\textheight,width=\textwidth,trim=.1cm .1cm .3cm .1cm, clip=true]{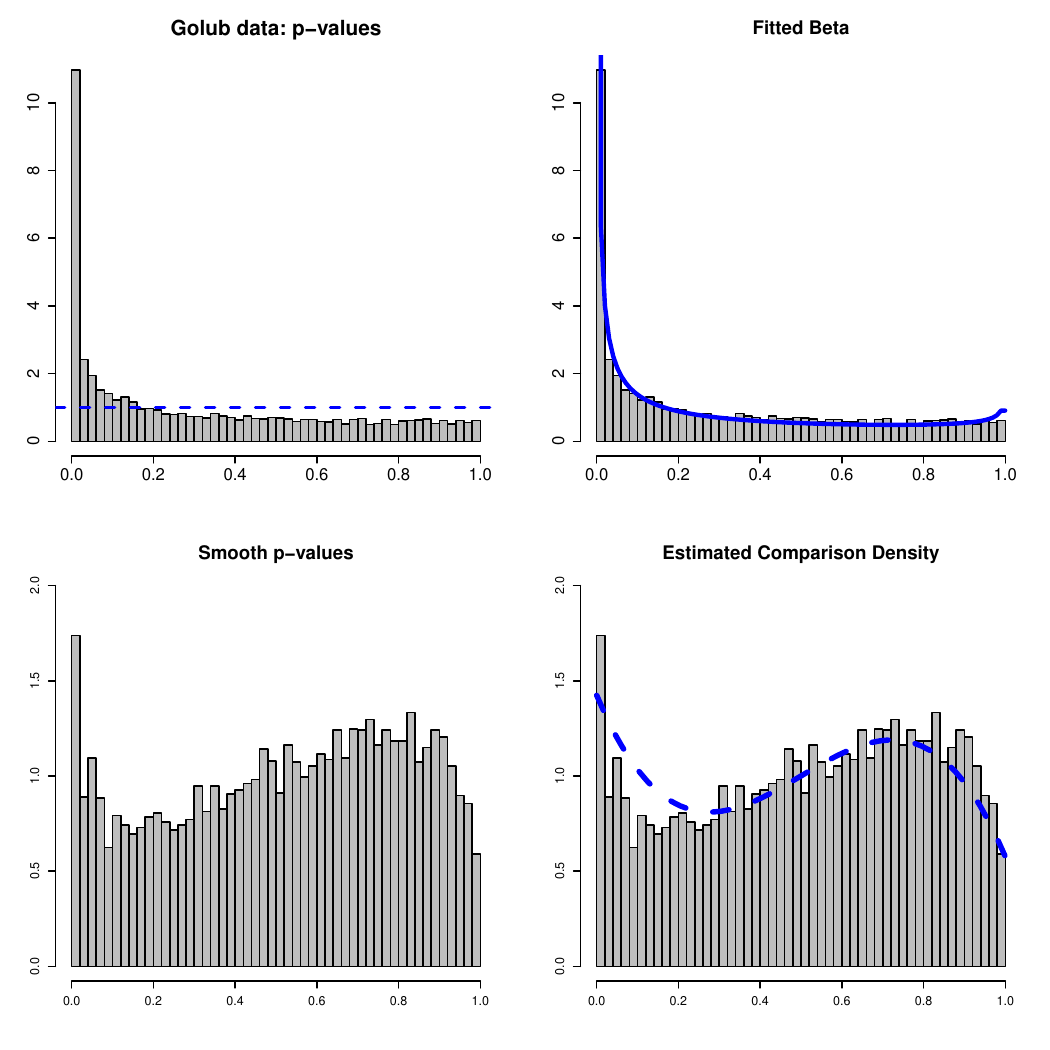}
\caption{The concepts of `smooth' p-value and the mechanism of skew-Beta density estimation for Golub gene expression data is illustrated; (A) [Top Left] Histogram of the $7129$ p-values of using two sample t-test; (B) [Top Right] Fitted beta distribution $\rm{Beta}(\hat \al=.32,\, \hat \be=.75)$ to the p-values; (C) [Bottom Left]  Distribution of the smooth p-values $v\,=\,F_{\rm{B}}(u;\,\hat \al=.32, \hat \be=.75)$; (D) [Bottom Right] LP-adaptive orthogonal (Nonparametric) density estimator of the smooth Golub p-values: $\dhat(v;\, F_{\rm{B}},F ) = 1 - .16 \Leg_3(v), 0<v<1$.}
\label{fig:golub}
\end{figure*}

Fig \ref{fig:golub}(A) shows the two-sample t-test p-values of $N=7129$ gene expressions of Golub cancer data \citep{golub99} while comparing $n_1=27$  acute lymphoblastic leukemia (ALL) and $n_2=11$ acute myeloid leukemia (AML) tumor samples to identify differentially expressed genes (data available in R package \texttt{golubEsets}). The comparison density, which is the distribution of p-values from this large-scale study indicates two nonparametric modeling challenges: First, we need density estimators for compact support $[0,1]$. Second, the sharp narrow peak near the boundary $0$ necessitates modeling of the highly dynamic tail.


Developing nonparametric density estimation for compact support $[0,1]$ is known to be a challenging problem due to the ``boundary effect.'' Recently, there has been a great deal of interest in adapting kernel density estimator for unit interval (see \cite{geenens2014} and \cite{wen2014}) by transforming the variable of interest into another one whose density estimation should be free from boundary problems, and then finally transforming that estimate back into the initial scale. However, it has been recognized that these methods are computationally not efficient. Other approaches, like regression-based density estimators via smoothing splines or local polynomials, are known to have a larger variance near the boundaries \citep{Thas10}.

For large-scale signal detection problems, besides compact support, the other more challenging modeling aspect is to tackle the highly dynamic tail near the boundaries $0$ and $1$, where all previously mentioned density estimators perform poorly.  We will propose a specially tailored \emph{new genre of density estimator} to tackle this typical non-standard ``U'' \emph{shaped density on the unit interval}.

We do note, however, that, it is not difficult to propose new density estimation techniques to fit the data, such as Fig \ref{fig:golub}(A), but it is less trivial to come up with a sparse parametric model that fits the data well and is easy to interpret. The most stunning fact about our approach (that we will elaborate on the next section) is that it require only \textit{three parameters} to accurately model the golub p-values, including the tail region!

To model this typical shape of comparison density (that arises in the large scale signal detection problems), we prefer estimators that are simultaneously accurate, computationally simple and parsimonious. Our proposal consists of two main steps: (1) converting ``spiky'' p-values to ``smooth'' p-values via the preflattening technique as described in Section 3.2 and (2) estimating  ``smooth'' p-values by expanding preflattened $d(u)$ (nonparametric $L^2$ or orthogonal series density estimator) or $\log d(u)$ (maximum entropy exponential density estimator) as a linear combination of orthonormal shifted Legendre polynomials (whose support is $[0,1]$), described in Section 3.3. \emph{The novelty of our approach lies in its unique ability to ``decouple'' the density estimation problem into two separate modeling problems: the tail part and the central part of the distribution.}
\subsection{Skew-Beta Density Estimator: Tail Modeling}
We propose a new genre of nonparametric density estimation technique, which starts with a parametric model $g(x)$ and permits the following \emph{universal decomposition}:
\beq
f(x) \,=\,g(x)\, \times\, d\big[ G(x); G, F\big].
\eeq
Verify comparison density \eqref{eq:cd} evaluated at $G(x)$ has the form $d(G(x);G,F)=f(x)/g(x)$. We call this new class as Skew-G density model for $f(x)$.

To efficiently capture the typical shape of the p-values shown in Fig \ref{fig:golub}(A), we select $g(x)$ on $[0,1]$ in such a way that it can tackle the rapidly changing tail. Hence, beta distribution is a good candidate for $g(x)$, which can act as pre-flattening function, as shown in Fig \ref{fig:golub}(B,C). The resulting skew-beta comparison density estimator has the following decomposition:
\beq \label{eq:beta}
d(u;F_0,F)\,=\,f_{\rm{B}}(u;\,\al,\be)\, \times\,d\big[ F_{\rm{B}}(u;\al,\be);\, F_{\rm{B}},F     \big], \quad {\rm for}~ 0<u<1,
\eeq
where beta density and distribution is denoted by $f_{\rm{B}}$ and $F_{\rm{B}}$. We define the quantity $F_{\rm{B}}(u;\al,\be)$ as the ``smooth'' p-value, shown in Fig \ref{fig:golub}(C), whose density is $d[ F_{\rm{B}}(u;\al,\be);\, F_{\rm{B}},F],$ $0<u<1$. Use simple maximum-likelihood estimate $\widehat \al$ and $\widehat \be$ to get the fitted beta; see Web Appendix E for more discussion and simple one line R implementation. Note that if none of the genes were differentially expressed for Golub data then we would expect the distribution of the pvalues to be $\text{Beta}(\al=1,\be=1)$, i.e., uniform. This simple observation can lead to a quick diagnostic for testing the presence of signals. Web Appendix F describes the procedure.

\subsection{Estimating density of ``Smooth'' P-values}
In this Section, we develop a nonparametric density estimator for smooth (Beta-transformed) p-values $d(v;F_{\rm{B}},F)$ on the unit interval, where $v=F_{\rm{B}}(u;\al,\be)$. Our proposal is based on $L^2$ or maximum entropy exponential comparison density estimator by expanding it in an orthonormal basis. We select shifted orthonormal \underline{L}egendre \underline{P}olynomials (LP) $\Leg_j(v)\,,$  $j=1,2,\ldots$ as basis, which form a complete orthonormal basis for $L^2[0,1]$ Hilbert space of square-integrable functions.
\bea
\text{$L^2$ orthogonal series estimator:} &  d(v;F_{\rm{B}},F) - 1 ~= ~\sum_{j} \LP[j;F_{\rm{B}},F] \Leg_j(v) \\
\text{Maximum entropy estimator:}& \log d(v;F_{\rm{B}},F)~= ~\sum_{j} \te_j \Leg_j(v) - K(\boldsymbol \te)
\eea
Due to simplicity and computational ease (which will be clear soon) we use $L^2$ estimate, which has nice asymptotic theoretical properties  \citep{ander80}. The large-$N$ paradigm (for golub data $N=7129$) makes all of these asymptotic analyses very much relevant and directly applicable. However, the reader is free to use the nonparametric maximum entropy exponential estimate (3.4).

\emph{Computation of LP-Fourier Coefficients}. The $L^2$ orthogonal expansion coefficients satisfies
\beq
\LP[j;F_{\rm{B}},F]\,=\,\int_0^1 d(v;F_{\rm{B}},F)\,\times\,\Leg_j(v) \dd v.
\eeq
Substitute $v$ by $F_{\rm{B}}(u;\al,\be)$ and note
\[d\big[ F_{\rm{B}}(u;\al,\be)\,   ;F_{\rm{B}},F\big]\,=\,\dfrac{d(u;F_0,F)}{f_B(u;\al,\be)}.\]
This lead to the following representation:
\beq \label{eq:LPcoef}
\LP[j;F_{\rm{B}},F]\,=\,\int_0^1 d(u;F_0,F)\, \Leg_j\big\{ F_{\rm{B}}(u;\al,\be) \big\} \dd u. \eeq
Eq. \eqref{eq:LPcoef} implies that the LP-coefficients can be represented as the mean of $\Leg_j$ \emph{evaluated at} the $F_{\rm{B}}(u;\al,\be)$ (these are Beta-transformed smooth pvalues). This representation result is summarized in the following theorem.
\begin{thm} The Fourier-LP coefficients admit the following representation
\[\LP[j;F_{\rm{B}},F]~=~\int_0^1 \Leg_j(v) \dd D(v;F_{\rm{B}},F)~=~\Ex\big[ \Leg_j\{F_{\rm{B}}(U;\al,\be)\}\big].\]
\end{thm}

As a practical consequence, this result allows fast computation of the LP expansion coefficients as mean of the $\Leg_j$ score functions \emph{evaluated at the transformed smooth pvalues}:
\beq
\LP[j;F_{\rm{B}},\widetilde F]\,\,~  \larrow ~\, \, N^{-1}\sum_{i=1}^N \Leg_j\big[ F_{\rm{B}}(u_i;\, \al,  \be) \big].
\eeq

\emph{Adaptive Estimation}. Estimate the smooth nonparametric model by selecting the `significantly large' LP-coefficients in a data-driven way. We will use the \cite{ledwina94} scheme using Schwarz selection criterion, which is known to be consistent. To identify the important coefficients using Ledwina's data-driven model selection criterion, we first rank the squared LP-coefficients. Then we take the penalized cumulative sum of $k$ coefficients using information criterion $N^{-1}\log(N)k$ and choose the $k$ for which this is the maximum. Fig \ref{fig:golub}(D) shows the resulting smooth estimate for the golub data.
\subsection{Estimating Proportion of True Null Hypothesis}
In Section 2.2 we have seen the thresholds for CD-based false discovery methods depend on the parameter $\pi_0$, the true proportion of noise or null hypothesis. The following is our proposed data-analytic estimation algorithm.
\vskip.3em
\textbf{Algorithm 1} [$\pi_0$ estimation by Minimum Deviance Criteria (MDC)]
\vskip.2em
\textit{Step 1.} Define $\cU_{\la}=\{u_i: \dhat(u;F_0,F) < \la\}$, where $\widehat{d}$ is the estimated beta-preflattened nonparametric comparison density; $|\cU_{\la}|=N_{\la}$. For each fixed $\la$ perform the following steps.
\begin{itemize}
 \item[(a)] Compute $\widetilde \LP_{\la}[j] \,\larrow \, N_\la^{-1} \sum_{i=1}^{N_\la} \Leg_j(u_i)$, which is the score coefficient for the following $L_2$ comparison density $\tilde d_{\la}(u)\,=\, 1\,+\, \sum_{j=1}^M \widetilde \LP_{\la}[j]\, \Leg_j(u)$,  based on $\cU_{\la}$.
\item[(b)] Calculate the deviance statistic $I_{\la} \,\larrow \, \sum_{j=1}^M \big| \widetilde \LP_{\la}[j]\big|^2$.
\end{itemize}
\textit{Step 2.} Display the deviance path $(\la, I_{\la})$ on a fine grid of $\la\in[1,\ga]$ and set $\la^{*} \,\larrow \, \arg\min_{\la} I_{\la}$.
\vskip.2em
\textit{Step 3.} Output $\widehat \pi_0 = N_{\la^{*}}/N$.
\vskip.35em

The \emph{rationale} behind the algorithm comes from the simple fact that under $H_0$, when all the cases are null, the underlying comparison density $d(u;F_0,F)$ should not deviate much from $\rm{Uniform}[0,1]$ as $D(u;F_0,F)=u$. Parseval's theorem dictates the following equality for $L^2$ comparison density:
\beq \label{eq:parseval} \int_0^1 \big[d(u;F_0,F) - 1\big]^2 \dd u\,\,=\,\,\sum_{j} \big| \LP[j;F_0,F] \,\big|^2. \eeq
In the light of Eq. \eqref{eq:parseval} the statistic $I_{\la}$ \emph{can be interpreted as the deviation of $\tilde d_{\la}(u)$ from uniformity}. Fig 1 of Web supplement illustrates this idea for Prostate cancer dataset (described in Section 4), where the shape of the estimated comparison density clearly indicates the presence of signals in the two tails. The deviance path for the rejection region of interest $1 \le \la \le \ga=3.5$ is shown in the right panel for $M=10$, which gives $\arg\min_{\la} I_{\la}=1.98$ and $\widehat \pi_0=0.971$. At the point $\la^{*}=1.98$ the deviance statistics takes the minimum value, which implies that the set $\cU_{\la=1.98}$ is most likely consists of the null cases.

\emph{Connection with other approaches}. Instead of working with comparison density, there are related methods based on comparison distribution function to identify p-values that differ significantly from the null $\rm{Uniform}[0, 1]$. Our method could be viewed as the formal density analogue of the graphical method proposed in \cite{schweder1982}. Another similar method proposed by \cite{storey10} estimates $\widehat \pi_0(\la) =(1-\widetilde D(\la))/(1-\la)$. They proposed a computationally intensive method to select the tuning parameter $\la$. It consists of bootstrapping p-values for each $\la$ in the range $\{0,.05,\ldots,.95\}$ and selecting the one which minimizes mean-squared error (MSE).
\subsection{Estimating Empirical Null}
Until this point we have assumed that the continuous null distribution $F_0$ is completely specified. As noted by \cite{efron2010book}, ``This is a less-than-usual occurrence'' for real datasets. \cite{efron04} proposed a method for estimating the unknown location and scale parameters under the zero assumption (i.e., the non-null or the mixing distribution (signal) $g(z)=0$ in a certain interval around the $z=0$, which is known to be prone to bias) for $F_0=\Phi[(x-\mu_0)/\si_0]$. \cite{Om10} estimates by putting a ${\rm Dirichlet}(\be)$ prior penalty on the $\pi$, where $\be$ is an extra tuning parameter. Here we propose a simple and fast solution that works quite nicely.

In the event all the observations $Z_1,\ldots,Z_N$ are coming from $F_0$, we have $F(z) = F_0[(z-\mu_0)/\si_0]$, where $F_0$ is the theoretical location-scale null distribution. We could have easily estimated the unknown $\mu_0$ and $\si_0$ just by fitting a simple linear regression model in the QQ plot. But instead, we have few $Z_i$'s, which are coming from $G$ ($\pi_0$-contamination neighborhood of $F_0$). If we knew which of these were unusual observations we could have removed them before fitting the linear regression. But as this is not the case, the natural step is to fit \emph{robust} linear regression instead of linear least square, which is resistant to the data contamination. To get the estimates, solve the M-estimating equation:
\bea
&\sum_{i=1}^n \Psi_T\Big[ \widetilde F^{-1}(u_i) - \mu_0 - \si_0 F_0^{-1}(u_i)\Big]~=~0,\nonumber \\
&\sum_{i=1}^n \Psi_T\Big[ \widetilde F^{-1}(u_i) - \mu_0 - \si_0 F_0^{-1}(u_i)\Big]F_0^{-1}(u_i)~=~0,
\eea
where we choose $\Psi_T$ to be the Tukey's biweight influence function \citep{Tukey1974} given by
\[ \Psi_T(z) = \left\{ \begin{array}{ll}
         z[ 1 - (z/k)^2]^2 & \mbox{for $|z| \le k$},\\
        0 & \mbox{for $|z|>k$}.\end{array} \right. \]
The value $k = 4.685$ is usually used as a default choice that provides 95\% asymptotic efficiency under normality and still offers protection against outliers. A straightforward R-implementation is possible via MASS library function \texttt{rlm}; see Web Appendix G for R-code. Table \ref{tab:enull} reports the findings of our approach for five large-scale studies discussed in Chapter 6 of \cite{efron2010book}. Finally, produce the null-adjusted nonparametric comparison density estimate using $\widehat F_0:=F_{0;\hat \mu_0, \hat \si_0}$
\beq
\widehat{d}(u;,\widehat F_0,F) \,=\, f_{\rm{B}}\big(\widehat F_0;\,\hat \al,\hat \be \big)\,\times\, \dhat \big( F_{\rm{B}}(\widehat F_0;\,\hat \al,\hat \be ) ;\, F_{\rm{B}},F   \big).
\eeq
\begin{table}
\caption{\textit{Comparing three empirical null estimation methods on five large scale data sets. Readers are referred to Appendix B of \cite{efron2010book} for more information on the data sets.}}
\vskip1em
\centering
\def\arraystretch{1.25}%
\begin{tabular}{ cccc }
\hline
\hline
Examples & Methods & $\widehat{\mu}_0$  & $\widehat{\si}_0$ \\ \hline
\multirow{3}{*}{HIV} & Robust biweight M-estimator & $.121$ & $.845$ \\
                     & MLE (Locfdr) & $.115$ & $.753$ \\
                     & Penalized mixture model (Mixfdr) & $.131$ & $.838$ \\ \hline
\multirow{3}{*}{Gene-Tagging} & Robust biweight M-estimator & $.383$ & $1.25$ \\
                     & MLE (Locfdr) & $.310$ & $1.26$ \\
                     & Penalized mixture model (Mixfdr) & $.373$ & $1.23$ \\ \hline
\multirow{3}{*}{Leukemia} &Robust biweight M-estimator& $.017$ & $1.88$ \\
                     & MLE (Locfdr) & $.120$ & $1.59$ \\
                     & Penalized mixture model (Mixfdr) & $.017$ & $1.86$ \\ \hline
\multirow{3}{*}{Prostate} & Robust biweight M-estimator & $-0.001$ & $1.092$ \\
                     & MLE (Locfdr) & $-0.002$ & $1.087$ \\
                     & Penalized mixture model (Mixfdr) & $-0.002$ & $1.068$ \\ \hline
\multirow{3}{*}{NYC Police} & Robust biweight M-estimator& $.087$ & $1.42$ \\
                     & MLE (Locfdr) & $.121$ & $1.39$ \\
                     & Penalized mixture model (Mixfdr) & $.080$ & $1.42$ \\ \hline

\hline
\end{tabular}
\vspace{-2em}
\label{tab:enull}
\end{table}

\section{Examples}

\subsection{Real Data Application}
Prostate cancer data \citep{prostate} consists of $102$ patient samples ($50$ labeled as normal and $52$ as prostate tumor samples) and $6033$ gene expression measurements. We aim to detect interesting genes that are differentially expressed in the two samples. For this purpose, we compute the two-sample t-test statistic $t_i$ for each gene and convert them into z-scale by $z_i~\leftarrow ~ \Phi^{-1}(\mathcal{T}_{100}(t_i))$, where $\mathcal{T}$ denotes the t-distribution function, shown in Fig 2(A) of Web Supplement.
Fig 2(D) of Web Supplement shows the final beta-preflattened smooth estimate of comparison density given by
\beq
\dhat(u;\Phi,F)\,=\, .68 \,\big[ 1\,+\, 0.057 \Leg_6 \big( F_{\rm{B}}(u;\widehat \al=.81,\widehat \be=.82) \big)  \big]\,u^{-.19}\,(1-u)^{-.18},~ ~~0<u<1,
\eeq
which along with the Minimum Deviance Criteria gives $\widehat \pi_0=.971$. The representation result \eqref{eq:locfdr} given in Theorem 3 immediately provides a CD-based estimate of  local fdr (which we call \texttt{CDfdr}). We compare our result with Locfdr and Mixfdr \citep{Om10} that estimate the fdr by separately estimating the numerator $\hf_0$ and the denominator $\hf$ (see Web Appendix D) . Naturally there are many variants available for these two methods depending on the way they estimate null and marginal densities. We have used the R package \texttt{locfdr} and \texttt{mixfdr} for implementation purpose. Locfdr estimates pool destiny $f$ using splines. Methods for estimating null include the following: (a) theoretical ($\cN(0,1)$); (b) maximum likelihood (MLE); (c) central matching (CM); (d) split-normal (SN). Mixfdr implement $J$ group normal mixture model for $f$. Estimation of empirical null involves putting Dirichlet prior on mixing proportion. We have used the default choice of $J$ and Dirichlet parameter $P$ throughout. All the three empirical null methods, including ours (reported in Table 1) shows that the theoretical null is not drastically different from $\cN(0,1)$. A very slight scale correction was required.
\vskip.25em
\begin{figure*}[!thb]
 \centering
 \includegraphics[height=.45\textheight,width=\textwidth,keepaspectratio,trim=.5cm 1cm .5cm .5cm, clip=true]{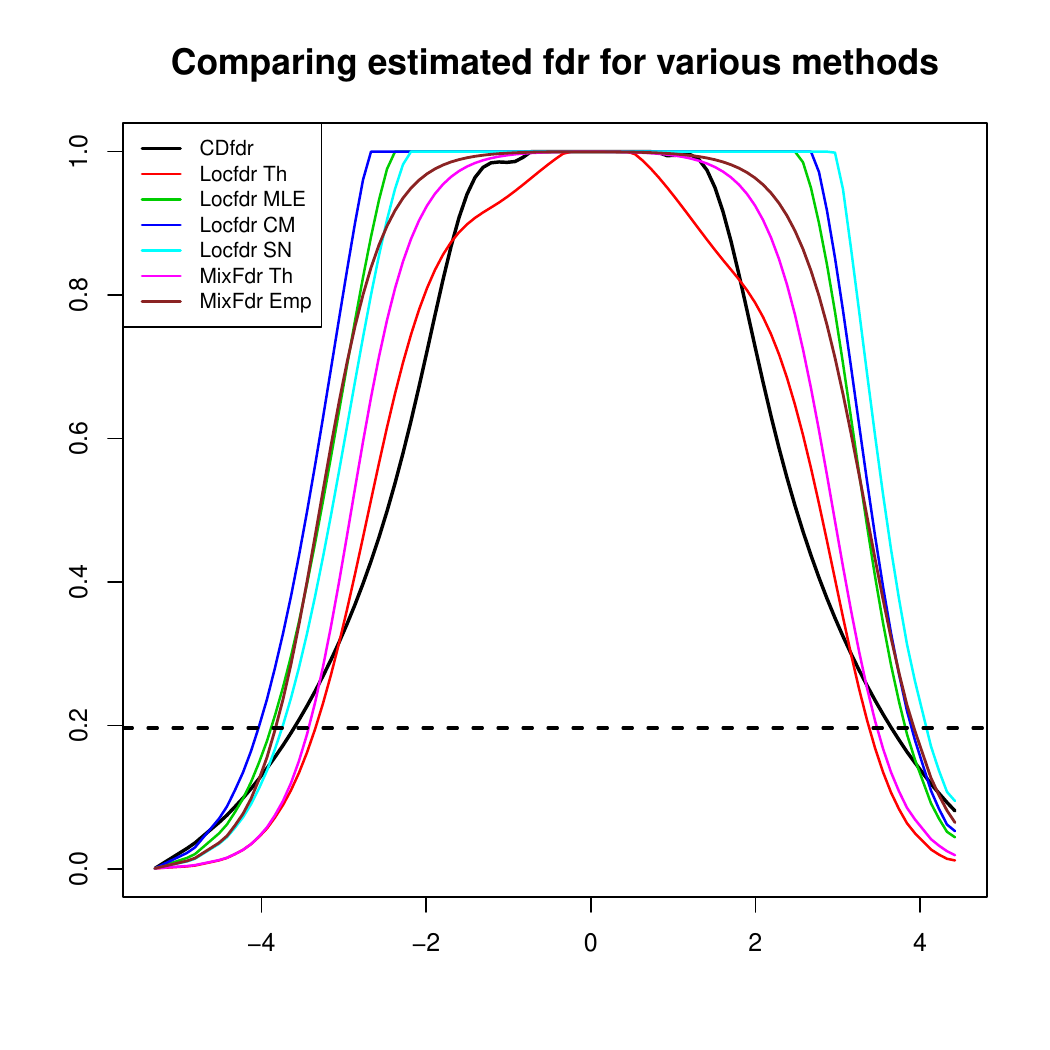} \\
\vspace{-.5em}
\caption{(color online) Estimated local fdr is shown for seven different methods for Prostate cancer data.}
\end{figure*}

Each the empirical null based estimated fdr functions, including ours, shows strong similarity in the tail-regions. The behavior of the fdr function in the tail-region (compared to the central region) is most crucial, as this determines which p-values will be declared as false or non-null. So, it might be more appropriate to focus on the quality of ``tail-modeling'' instead of considering the entire curve. We carry out this in the next section where we carefully quantify the estimation accuracy, especially in the tails.

The number of non-null genes identified by different methods, along with the estimates of the proportion of true null is provided in the Web Supplement Table 2. First note that the estimates of $\pi_0$ from Locfdr-CM and Locfdr-SN are unrealistic as they exceed $1$. Number of non-null genes identified using CDfdr matches with the Locfdr-SN and Mixfdr-Emp method. Furthermore, our proposed minimum deviance estimate $\widehat{\pi}_0$ is very close to the Mixfdr-Th, Mixfdr-Emp and Locfdr-Th methods.
\subsection{Simulation Study}
In order to further evaluate the accuracy of the \texttt{CDfdr} algorithm, we perform two simulated experiments. We are mainly interested in investigating how accurately different methods estimate fdr, especially in the tail-regions based on the mean integrated square error (MISE) criteria. Comparisons will be done with Locfdr, Mixfdr and Fdrtool \citep{strimmer}. Grenander density estimation is used in Fdrtool for  estimating the unconditional density $f$ and is implemented in R package \texttt{Fdrtool}.
\subsubsection{Mixture Normal}
We simulate $T_i \sim \cN(\mu_i,1),i=1,\ldots,N=5000$ out of which $4500$ $\mu_i$ is set to zero. The remaining $500$ is drawn from (once and for all) $\cN(\mu,1)$. We estimate the $\widehat \fdr$ for various methods and repeat the whole process $150$ times for $\mu=0.2,0.5,1,2$. Our setup closely follows \cite{storey10} and \cite{Om10}.

\emph{Local fdr estimation}. The goal is to investigate the efficiency of different estimation methods when we fix the null density at $\cN(0,1)$.   Web Supplement Figure 3 depicts the expectation and standard deviation of local fdr for various methods under four different choices of $\mu$. For $\mu=2$, when the signal and noise are well-separated, all of the methods perform equally well apart from Fdrtool, which not only shows high bias but has large variability in the crucial tail region. Under the more difficult scenario of $\mu=.2$, clearly CDfdr is the only method that can claim to be unbiased. The variability of Mixfdr and CDfdr seems similar, though in the extreme (right) tail, CDfdr shows more stability. If we look at the $\rm{Sd}(\widehat \fdr)$ curves for CDfdr and compare with other competing methods, it appears to be the least variable. Also, CDfdr achieves \emph{near unbiasedness irrespective of the underlying signal strength}, which makes it a reliable tool for large-scale inference problems.

\emph{Null proportion estimation}. It is interesting to examine how all of these methods would perform in estimating the true null proportion under different degrees of signal strength. We implemented our Minimum Deviance Criteria (MDC) algorithm to estimate $\pi_0$. The simulation was done under the same experimental setup. The boxplots of estimates of $\pi_0$ under different non-null densities are shown in Web Supplement Figure 4. Locfdr shows the largest variability. Mixfdr certainly performs best among the competing methods. However, the method that was particularly successful for estimating the true null proportion $\pi_0=.9$ \emph{quickly and accurately} is based on the our proposed MDC. This makes the CDfdr algorithm more powerful and efficient as a signal detection tool that protects against false discovery.
\subsubsection{Mixture Uniform}
We generate p-values from the following model with the parameter choices: $\pi_0=\{0.9,0.95,0.99\}$ and $a=\{0.02,0.002\}$:
\beq
\pi_0\, \rm{Uniform}[0,1]\,+\, (1-\pi_0) \, \rm{Uniform}[0,a].
\eeq
Here we particularly pay attention to the tails and for that we consider the tail-specific MISE criteria $\E \int_{\mathcal{S}} (\fdr(u) - \hat \fdr(u))^2 \dd u$, where $\mathcal{S}$ denotes the collection of $u$ coming from the alternative model $U[0,a]$. The goal is to quantify how precisely the fdr is estimated for the signals (in the tail). Here the parameter $a$ controls the signal strengths and parameter $\pi_0$ determines the underlying sparsity levels.
\begin{figure}[thb]
\vspace{-1em}
 \centering
 \includegraphics[height=\textheight,width=\textwidth,keepaspectratio]{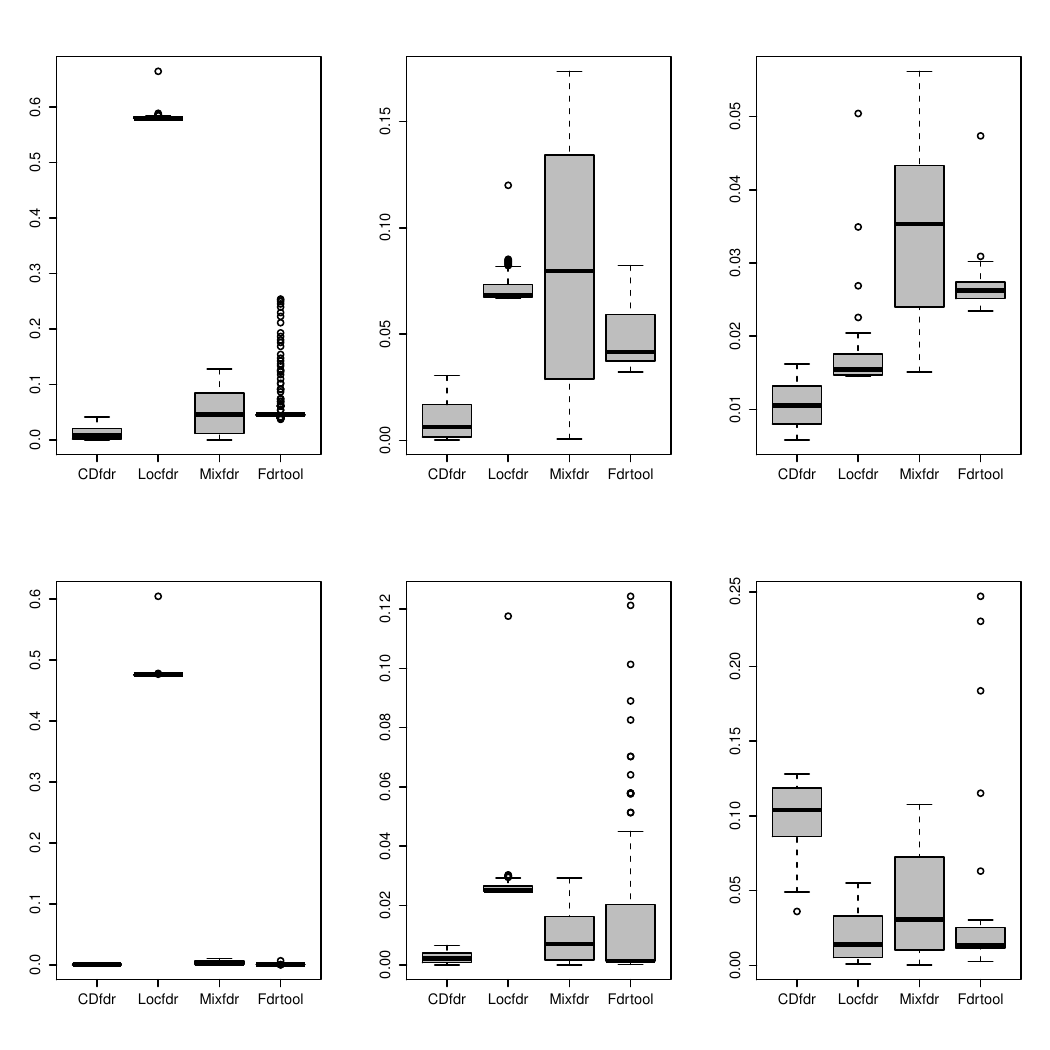}
\vspace{-1.5em}
\caption{Compares the tail specific MISE for Uniform mixture model $\pi_0\, \rm{Uniform}[0,1]\,+\, (1-\pi_0) \, \rm{Uniform}[0,a]$. The rows corresponds to $a=0.02$ and $0.002$. For each row the columns (from left to right) denotes $\pi_0=0.9,0.95,0.99$.}
\label{fig:simu:unif}
\end{figure}
Our simulation design covers the complete spectrum from
\[ \text{dense and weak}\, \rarrow \, \text{rare and weak} \, \rarrow  \, \text{strong and dense} \, \rarrow \, \text{strong and sparse signal}.  \]
In the presence of weak signals ($a=0.02$), the Locfdr and Mixfdr show a great deal of variability, as shown in Fig \ref{fig:simu:unif} (first row). CDfdr maintains the smallest tail-specific MISE among all the methods, which again ensures its utility.  For strong signals (second row of Fig \ref{fig:simu:unif}), most of the methods have reasonable performance, except perhaps, $(a=0.002,\pi_0=0.9)$ case where the Locfdr poorly approximates the tail. The large number of outliers for Fdrtool are also not desirable.
\vskip.25em
Overall, it is encouraging that \texttt{CDfdr} \emph{adapts to the underlying signal sparsity and strength in many cases}, which makes it very attractive and reliable for large-scale studies. Undoubtedly for the examples we have discussed in this paper it appears that CDfdr shows the most consistent and robust performance.
\section{Conclusion}
We have shown how the concepts and notations of comparison density and distribution (1) provide unification (and an alternative rationale) of the two cultures of multiple testing: Frequentist BH procedure and Efron's empirical Bayes local fdr idea; and (2) lead to a new class of efficient nonparametric signal detection algorithms. Our approach converts the mixed sample signal detection problem into a nonparametric comparison density function estimation problem. The author expects proposed CD-based nonparametric functional statistical approach will not only simplify the theory (which has an enormous literature) and practice but will also allow the possibility to include large-scale inference topics \emph{as part of a conventional statistical inference curriculum} that can cover one-sample, mixed-sample and two-sample problems using \emph{single concept, notion, and algorithm}. Our unified treatment by \emph{linking} different large-scale multiple testing ideas (rather than presenting them as ``unrelated tools'') could accelerate the learning for students.

We made a crucial observed that almost all large-scale inference problems produce comparison density with a typical stylized shape -- ``U'' shaped density over the compact interval $[0,1]$, which is not easily estimable using traditional nonparametric density estimation techniques. To address this, we introduced a \emph{new genre} of density estimation methods based on the idea of Skew-G density decomposition. This allows for richer data-driven specification for tail-modeling via simple parametric models, which has an added advantage of being interpretable and easily implementable. We believe this technique can be used in many other modeling problems (outside multiple testing), such as heavy-tailed density estimation.

In future work, we would like to extend this technique to discrete mixed sample problems. Although we believe the theory presented in \cite{Deep14LP} will guide us in that direction, it will require more detailed studies. Systematic theoretical investigation of the asymptotic properties of the proposed nonparametric skew-beta density estimator is an unexplored and open problem. Additionally, we plan to understand the effect of dependence on our proposed density estimator and how it influences the false discovery thresholds.
\section{Supplementary Materials}
Web Appendices, Tables and Figures referenced in Sections 1.1, 2.2, 3.2, 3.4 and 4 are available with this paper at the \emph{Biometrics} website on Wiley Online Library. Online link: {\color{blue} http://bit.do/LSSD-BiometricsSupp}.

\section*{Acknowledgment}
The author would like to thank Emanuel Parzen for several valuable comments and suggestions. The author also thank the Associate Editor and anonymous referees whose constructive comments have greatly helped to improve the quality and presentation of the paper. The author dedicates this paper to John Tukey, the pioneer of multiple comparison idea, on the occasion of his 100th birthday.

\bib

\end{document}